\documentclass[12pt,amssymb,verbatim]{amsart}

\usepackage[all]{xy}

\def\cal{\mathcal}

\numberwithin{equation}{section}

\theoremstyle{plain}
\newtheorem{theorem}{Theorem}

\newtheorem{proposition}[theorem]{Proposition}
\newtheorem{lemma}[theorem]{Lemma}
\newtheorem{corol}[theorem]{Corollary}

\theoremstyle{definition}

\theoremstyle{remark}
\newtheorem{remark}[theorem]{Remark}

\newtheorem{example}[theorem]{Example}

 \def\today{\ifcase\month\or
  January\or February\or March\or April\or May\or June\or
  July\or August\or September\or October\or November\or December\fi
  \space\number\day, \number\year}

\begin{document}
\title[Roots of polynomials]
{Roots of polynomials over division
rings}
\author[Goutor, Tikhonov]{Alina G. Goutor, Sergey V. Tikhonov}


\address{Goutor:
Belarusian State University, Nezavisimosti Ave., 4,
220030, Minsk, Belarus} \email{goutor@bsu.by }
\address{Tikhonov:
Belarusian State University, Nezavisimosti Ave., 4,
220030, Minsk, Belarus} \email{tikhonovsv@bsu.by }

\def\cA{{\cal A}}
\def\cB{{\cal B}}
\def\cC{{\cal C}}
\def\cD{{\cal D}}
\def\cV{{\cal V}}
\def\cE{{\cal E}}
\def\cM{{\cal M}}
\def\cR{{\cal R}}
\def\cQ{{\cal Q}}
\def\M{{\rm M}}
\def\cS{{\cal S}}
\def\Symb{\textrm{Symb}}
\def\Gal{\textrm{Gal}}
\def\U{\textrm{U}}
\def\SU{\textrm{SU}}
\def\R{\textrm{R}}
\def\GL{\textrm{GL}}

\def\U{{\rm U}}
\def\SU{{\rm SU}}
\def\SL{{\rm SL}}
\def\op{{\rm op\;}}


\def\cor{\textrm{cor}}
\def\deg{\textrm{deg}}
\def\exp{\textrm{exp}}
\def\Gal{\textrm{Gal}}
\def\ram{\textrm{ram}}
\def\Spec{\textrm{Spec}}
\def\Proj{\textrm{Proj}}
\def\Perm{\textrm{Perm}}
\def\coker{\textrm{coker\,}}
\def\Hom{\textrm{Hom}}
\def\im{\textrm{im\,}}
\def\ind{\textrm{ind}}
\def\int{\textrm{int}}
\def\inv{\textrm{inv}}
\def\min{\textrm{min}}

\def\w{\widehat}

\begin{abstract}
In this paper, we study properties of polynomials over division rings.
Moreover, we present formulas for finding roots
of some polynomials.
\end{abstract}

\maketitle


\def\dd{{\partial}}

\def\into{{\hookrightarrow}}

\def\emptyset{{\varnothing}}

\def\alp{{\alpha}}  \def\bet{{\beta}} \def\gam{{\gamma}}
 \def\del{{\delta}}
\def\eps{{\varepsilon}}
\def\kap{{\kappa}}                   \def\Chi{\text{X}}
\def\lam{{\lambda}}
 \def\sig{{\sigma}}  \def\vphi{{\varphi}} \def\om{{\omega}}
\def\Gam{{\Gamma}}  \def\Del{{\Delta}}  \def\Sig{{\Sigma}}
\def\ups{{\upsilon}}


\def\A{{\mathbb A}}
\def\F{{\mathbb F}}
\def\Q{{\mathbb{Q}}}
\def\CC{{\mathbb{C}}}
\def\PP{{\mathbb P}}
\def\R{{\mathbb R}}
\def\Z{{\mathbb Z}}
\def\X{{\mathbb X}}

\def\Gm{{{\Bbb G}_m}}
\def\Gmk{{{\Bbb G}_{m,k}}}
\def\GmL{{\Bbb G_{{\rm m},L}}}
\def\Ga{{{\Bbb G}_a}}

\def\Fb{{\overline F}}
\def\Hb{{\overline H}}
\def\Kb{{\overline K}}
\def\Lb{{\overline L}}
\def\Yb{{\overline Y}}
\def\Xb{{\overline X}}
\def\Tb{{\overline T}}
\def\Bb{{\overline B}}
\def\Gb{{\overline G}}
\def\Vb{{\overline V}}

\def\kb{{\bar k}}
\def\xb{{\bar x}}

\def\Th{{\hat T}}
\def\Bh{{\hat B}}
\def\Gh{{\hat G}}

\def\Xt{{\tilde X}}
\def\Gt{{\tilde G}}

\def\gg{{\mathfrak g}}
\def\gm{{\mathfrak m}}
\def\gp{{\mathfrak p}}
\def\gq{{\mathfrak q}}

\def\textrm#1{\text{\rm #1}}

\def\res{\textrm{res}}
\def\cor{\textrm{cor}}
\def\char{\textrm{char}}
\def\R{\textrm{R}}

\def\tors{_{\textrm{tors}}}      \def\tor{^{\textrm{tor}}}
\def\red{^{\textrm{red}}}         \def\nt{^{\textrm{ssu}}}
\def\sc{^{\textrm{sc}}}
\def\sss{^{\textrm{ss}}}          \def\uu{^{\textrm{u}}}
\def\ad{^{\textrm{ad}}}           \def\mm{^{\textrm{m}}}
\def\tm{^\times}                  \def\mult{^{\textrm{mult}}}
\def\tt{^{\textrm{t}}}
\def\uss{^{\textrm{ssu}}}         \def\ssu{^{\textrm{ssu}}}
\def\cf{^{\textrm{cf}}}
\def\ab{_{\textrm{ab}}}

\def\et{_{\textrm{\'et}}}
\def\nr{_{\textrm{nr}}}

\def\SB{\textrm{SB}}

\def\<{\langle}
\def\>{\rangle}

\def\til{\;\widetilde{}\;}


\font\cyr=wncyr10 scaled \magstep1%
\def\Bcyr{\text{\cyr B}}
\def\Sh{\text{\cyr Sh}}
\def\Ch{\text{\cyr Ch}}

\def\lpsi{{{}_\psi}}
\def\bks{{\backslash}}

\def\Br{\textrm{Br}}
\def\Pic{\textrm{Pic}}
\def\Bt{{{}_2\textrm{Br}}}
\def\Bn{{{}_n\textrm{Br}}}
\def\Grass{\textrm{Grass}}



\section{Introduction and preliminary results} \label{sec:intro}

In this paper, we study polynomials over division rings.
Let $\cD$ be an associative division ring. Let also $\cD[x]$ denote the polynomial ring in one variable $x$ over $\cD$, where $x$ commutes elementwise with $\cD$. The coefficients of such polynomials may not commute with elements of the ring.
Polynomials in $\cD[x]$ are added in the obvious way, and multiplied
according to the rule
$$
(a_nx^n + \dots +a_0)(b_mx^m + \dots +b_0)=(c_{m+n}x^{m+n} + \dots +c_0),
$$
where $c_k=\sum_{i+j=k}a_ib_j$.
For references on polynomial rings over division rings, see \cite[Ch. 5, \S 16]{La91} and \cite{GoMo65}.
The degree of $P(x)\in \cD[x]$ is defined in the usual way.
For a polynomial
$$
P(x)=a_nx^n + a_{n-1}x^{n-1} + \dots + a_1x+a_0 \in \cD[x]
$$
and an element $a\in \cD$, we define $P(a)$ to be the element
$$
a_na^n + a_{n-1}a^{n-1} + \dots + a_1a+a_0.
$$

An element $a \in \cD$ is said to be a (right) root of $P(x)$ if $P(a)=0$.
The noncommutative form of the Remainder Theorem says that an element $a\in \cD$ is a root of a nonzero polynomial
$P(x)$ iff $x-a$ is a right divisor of $P(x)$ in $\cD[x]$ (see, e.g., \cite[Prop. 16.2]{La91}).

For $a\in \cD$, the set
$$
[a] := \{dad^{-1} | d \in  \cD \setminus \{0\} \}
$$
will be called the conjugacy class of $a$. The centralizer of $a$ is defined as
$$
Z(a) := \{b \in \cD | ab=ba\}.
$$

Over a field, a polynomial of degree $n$ has at most $n$ distinct roots.
Over a division ring this is no longer true, but
by Gordon-Motzkin theorem (\cite[Th. 16.4]{La91}), a polynomial of degree $n$ in $\cD[x]$ has roots in at most $n$ conjugacy classes of $\cD$,
moreover, if $P(x)=(x-a_1)\dots(x-a_n)$, where $a_1,\dots,a_n \in \cD$, then any root of $P(x)$ is conjugate to some $a_i$.
Note that from $P(x)=L(x)R(x) \in \cD[x]$ it does not follow that $P(a)=L(a)R(a)$ (see Proposition \ref{prop:P=LR} bellow).
In particular, if $a$ is a root of $L(x)$, then $a$ is not necessarily a root of $P(x)$.

The problem of finding the roots of a polynomial over a division ring has been investigated in ring theory and applied mathematics.
The most studied is the case of polynomials over Hamilton's quaternion algebra $\mathbb{H}$.
In analogy to field theory, the notion of a (right) algebraically closed division ring $\cR$ is defined.
This is equivalent to saying that every polynomial in $\cR[x]$ splits completely into a product of linear factors in $\cR[x]$.
By Niven-Jacobson theorem (\cite[Th. 16.14]{La91}), the quaternion division algebra over a real-closed field is algebraically closed.
Baer's theorem (\cite[Th. 16.15]{La91}) says that noncommutative centrally finite (right) algebraically closed division rings are precisely the division rings of quaternions over real-closed fields. In \cite{HuSo02}, a formula was found for roots of any quadratic polynomial in $\mathbb{H}[x]$.
This formula was generalized to any quaternion algebra in \cite{Ab09} and \cite{Ch15}.
In \cite{SePeVi01}, it was shown that the roots of any polynomial in $\mathbb{H}[x]$ are roots of the real companion polynomial.
In \cite{JaOp10}, it was presented an algorithm for finding all roots of a polynomial in $\mathbb{H}[x]$ using
the real companion polynomial. In \cite{ChMa16}, a few of these results were generalized to the case of any central division algebra.
Recall that a central division algebra is a division algebra which is finite dimensional over its center.  A complete method for finding the roots of all polynomials over an octonion division algebra was described in \cite{Ch20}.

In \cite{FaMiSeSo17}, it was presented the following explicit formula describing roots of a product of linear factors in $\mathbb{H}[x]$.

\begin{theorem} \label{th:roots_quaternions}(\cite[Th. 4]{FaMiSeSo17}).
Let $P(x) = (x-q_n)\dots(x-q_1)$, where $q_1,\dots,q_n \in \mathbb{H}$.  If the conjugacy classes $[q_k]$ are distinct, then the polynomial $P(x)$
has exactly $n$ roots $\zeta_k$ which are related to the elements $q_k$ as follows:
$$
\zeta_k = \overline{P}_k(q_k)q_k(\overline{P}_k(q_k))^{-1}; k=1,\dots,n,
$$
$$
P_k(x) := \left\{
\begin{array}{l}
1,  \mbox{ if } k=1,\\
(x-q_{k-1})\dots(x-q_1), \mbox{ otherwise } \\
\end{array}
\right.
$$
and $\overline{P}_k(x)$ is the conjugate polynomial of $P_k(x)$.
\end{theorem}


In theorem \ref{th:roots_division} below, we generalize this formula for the case of any division ring.
This is the first aim of this paper.

Let $F$ be the center of a division ring $\cD$.
If $a$ is a root of a polynomial $f(x)\in F[x]$,
then any element from the conjugacy class $[a]$ is a root of $f(x)$. The conjugacy class $A$ is called algebraic over $F$ if one (an hence all) of
its elements is algebraic over $F$. If $A$ is algebraic over $F$, then the minimal polynomial of $A$ is, by definition, the minimal polynomial of any element from $A$.

In the case of a quadratic minimal polynomial, there is the following

\begin{theorem} \cite[Lm. 16.17]{La91}.
Let $\cD$ be a division ring with center $F$, and let $A$ be a conjugacy class of $\cD$ which
has a quadratic minimal polynomial $\lambda(x)$ over $F$. If $P(x)\in \cD[x]$ has two roots in $A$, then
$P(x)\in \cD[x]\lambda(x)$ and $P(x)$ vanishes identically on $A$.
\end{theorem}

This means that a polynomial over a quaternion division algebra may have two different types of roots:
isolated and spherical roots.
A root $q$ of $P(x)$ is called spherical if $q$ is not central and for every $d \in [q]$ we have
$P(d) = 0$. 
A root $q$ is called isolated if the conjugacy class $[q]$ contains no other root of $P(x)$.

The second aim of this paper is to show that in the case of a conjugacy class with minimal polynomial of bigger degree the situation is completely different.
More precisely, we proved the following


\begin{theorem} \label{th:quadratic}
Let $\cD$ be a noncommutative division ring with the center $F$, $a\in \cD$ an algebraic over $F$ element with minimal polynomial $\lambda(x)$ of degree $n > 2$.
Then there exists a quadratic polynomial $P(x) \in \cD[x]$ such that

1. $P(x)$ has infinitely many roots in the conjugacy class $[a]$,

2. there are infinitely many elements in $[a]$ which are not roots of $P(x)$,    

3. $\lambda(x)$ does not divide $P(x)$.
\end{theorem}



\section{Proof of Theorem \ref{th:quadratic}} \label{sec:proof}

In the proof of Theorem \ref{th:quadratic} we will use the following statements.


\begin{theorem} \label{th:vanish_on_class} (\cite[Th. 16.6]{La91}).
Let $\cD$ be a division ring with center $F$ and $A$ a conjugacy class of $\cD$ which is algebraic over $F$ with minimal polynomial $f(x)\in F[x]$.
A polynomial $P(x) \in \cD[x]$ vanishes identically on $A$ iff $P(x) \in \cD[x]f(x)$.
\end{theorem}


\begin{theorem} \label{th:two roots} (\cite[Th. 16.11]{La91} and \cite[Th. 4]{GoMo65}).
If a polynomial $P(x)\in \cD[x]$ has two distinct roots in a conjugacy class of $\cD$, then it has infinitely many roots in that class.
\end{theorem}


\begin{proposition} \label{prop:P=LR}(\cite[Pr. 16.3]{La91}).
Let $\cD$ be a division ring and let $P(x)=L(x)R(x)\in \cD[x]$. Let $d \in \cD$ be such that $h:=R(d)\ne 0$.
Then
$$
P(d)=L(hdh^{-1})R(d).
$$
In particular, if $d$ is a root of $P(x)$ but not of $R(x)$, then $hdh^{-1}$ is a root of $L(x)$.
\end{proposition}


\noindent {\it {Proof of Theorem \ref{th:quadratic}}}.
Let $a\in \cD$ be an element with minimal polynomial $\lambda(x)$ of degree $n>2$.
Let also $d \in \cD$ be an element such that $d$ does not commute with $a$.  Let $q=dad^{-1}$ and
$b=(q-a)q(q-a)^{-1}$. Then $q\ne a$ and $q \in [a]$.  By Proposition \ref{prop:P=LR}, $q$ is a root of the polynomial
$P(x) := (x-b)(x-a)$.

Since $a$ is also a root of $P(x)$, then by Theorem \ref{th:two roots}, $P(x)$ has infinitely many roots in $[a]$.
Moreover, since the degree of $\lambda(x)$ is bigger than 2, then $\lambda(x)$ does not divide $P(x)$.
Hence by Theorem \ref{th:vanish_on_class}, $P(x)$ does not vanish identically on $[a]$.

Suppose that $tat^{-1} \in [a]$ is not a root of $P(x)$. This means that $tat^{-1} \ne a$ and
$$
(tat^{-1} - a) tat^{-1} (tat^{-1}- a)^{-1} \ne b
$$
by Proposition \ref{prop:P=LR}.

Note that
$$
(tat^{-1} - a) tat^{-1}  (tat^{-1}- a)^{-1} = (t a - a t) a (t a -a t)^{-1}.
$$

Let $z \in Z(a)$, $t_1=t+z$ and $q_1=t_1 a t_1^{-1}$. Then
$$
(q_1-a)q_1(q_1-a)^{-1} = (t_1 a t_1^{-1}-a)t_1 a t_1^{-1}(t_1 a t_1^{-1}-a)^{-1}=
(t_1 a -a t_1) a (t_1 a -a t_1)^{-1}=
$$
$$
=((t+z) a - a (t+z)) a ((t+z) a - a (t+z))^{-1}=
(ta-at)a(ta-at)^{-1} \ne b.
$$
Thus $q_1$ is not a root of $P(x)$ by Proposition \ref{prop:P=LR}.

Now let $z_1 \in Z(a)$, $z_1 \ne z$.
Note that the centralizer $Z(a)$ is infinite by \cite[Th. 3]{GoMo65}.
Let also $t_2=t+z_1$ and $q_2=t_2 a t_2^{-1}$. Then
$q_2$ is also not a root of $P(x)$. 

Assume that $q_1=q_2$. Then
$$
(t+z)a(t+z)^{-1} = (t+z_1)a(t+z_1)^{-1} \Leftrightarrow
$$
$$
(t+z_1)^{-1}(t+z)a = a(t+z_1)^{-1}(t+z) \Leftrightarrow
$$
$$
(t+z_1)^{-1}(t+z_1 +(z-z_1))a = a(t+z_1)^{-1}(t+z_1 +(z-z_1)) \Leftrightarrow
$$
$$
a+ (t+z_1)^{-1}(z-z_1) a = a + a(t+z_1)^{-1}(z-z_1) \Leftrightarrow
$$
$$
(t+z_1)^{-1}(z-z_1)a = a(t+z_1)^{-1}(z-z_1) \Leftrightarrow (t+z_1)^{-1}a(z-z_1) = a(t+z_1)^{-1}(z-z_1) \Leftrightarrow
$$
$$
 (t+z_1)^{-1}a = a(t+z_1)^{-1}\Leftrightarrow a(t+z_1) = (t+z_1) a \Leftrightarrow at = ta.
$$
This gives a contradiction since $t$ does not commute with $a$. Then $q_1\ne q_2$.
Hence any $z\in Z(a)$ defines the element $(t+z) a (t+z)^{-1} \in [a]$ which is not a root of $P(x)$ and
all such elements are distinct. Since the centralizer $Z(a)$ is infinite,
then there are infinitely many elements in $[a]$ which are not roots of $P(x)$.
\qed



\begin{remark}
In the notation of the proof of Theorem \ref{th:quadratic}, let
$$
b_1 = (tat^{-1} - a) tat^{-1} (tat^{-1}- a)^{-1}.
$$
The polynomials $(x-b)(x-a)$ and $(x-b_1)(x-a)$ have infinitely many roots in the conjugacy class $[a]$, but
$a$ is the unique common root of these polynomials.

\end{remark}

\section{Roots of polynomials} \label{sec:roots}


It seems to us that the following lemma may be a known result, but we have not found an exact reference.
For the reader's convenience, we provide a proof here.


\begin{lemma} \label{lm:linear factors}
Let $\cD$ be a division ring.
Let also
$$
P(x) = (x-d_n)\dots(x-d_1),
$$
where $d_1,\dots,d_n \in \cD$.
If the conjugacy classes $[d_k]$ are distinct, then the polynomial $P(x)$ has exactly $n$ roots
and any root of $P(x)$ is conjugate to some $d_i$.
\end{lemma}


\noindent {\it {Proof}}.
By Gordon-Motzkin theorem (\cite[Th. 16.4]{La91}), the roots of $P(x)$ lie in $n$ conjugacy classes of $\cD$ and
any root of $P(x)$ is conjugate to some $d_i$. Let $d\in \cD$ be a root of $P(x)$, then
$P(x)=L(x)(x-d)$ for some $L(x)\in \cD[x]$. By Proposition \ref{prop:P=LR}, all roots of $P(x)$ different from $d$
are conjugate to roots of $L(x)$.
Since the conjugacy classes $[d_k]$ are distinct and  $\deg (L(x)) = n-1$, then by Gordon-Motzkin theorem,
$L(x)$ has no roots in $[d]$. Thus $P(x)$ has only one root in each conjugacy class.
\qed


\begin{theorem} \label{th:roots_division}
Let $\cD$ be a division ring with center $F$.
Let also
$$
P(x) = (x-d_n)\dots(x-d_1),
$$
where $d_1,\dots,d_n \in \cD$. Assume that $d_1,\dots,d_{n-1}$ are algebraic over $F$.
Let also $f_i(x)$ be the minimal polynomial of $d_i$, $i=1,\dots,n-1$.
If the conjugacy classes $[d_k]$ are distinct, then the polynomial $P(x)$
has exactly $n$ zeros $\zeta_k$ which are related to the elements $d_k$ as follows:
$$
\zeta_k = P_k(d_k)d_k(P_k(d_k))^{-1}; k=1,\dots,n,
$$
$$
P_k(x) := \left\{
\begin{array}{l}
1,  \mbox{ if } k=1,\\
S_1(x) \dots S_{k-1}(x), \mbox{ otherwise, } \\
\end{array}
\right.
$$
where $S_i(x)\in \cD[x]$ is such that $f_i(x)=S_i(x)(x-d_i)$, $i=1,\dots,n-1$.
\end{theorem}


\noindent {\it {Proof}}. Since $S_i(x)$ has coefficients in the field $F(d_i)$,
then
$$
f_i(x)=S_i(x)(x-d_i)=(x-d_i)S_i(x)
$$
for $i=1,\dots,n-1$. Note that
$$
P(x)P_k(x) = (x-d_n)\dots(x-d_1) S_1(x)\dots  S_{k-1}(x) =
$$
$$
=(x-d_n)\dots(x-d_k) f_{k-1}(x)\dots f_1(x).
$$
Since $f_{k-1}(x)\dots f_1(x) \in F[x]$, then $d_k$ is a root of the polynomial $P(x)P_k(x)$.
Note that for $i=1,\dots,k-1$, $d_k \not \in [d_i]$, so $d_k$ is not a root of $f_i(x)$ by Dickson's Theorem
(\cite[Th. 16.8]{La91}). Hence $d_k$ is not a root of $f_{k-1}(x)\dots f_1(x)$.
Then $d_k$ is not a root of $P_k(x)$. Indeed,
$$
(x-d_{k-1})\dots(x-d_1) P_k(x) = (x-d_{k-1})\dots(x-d_1) S_1(x)\dots  S_{k-1}(x) = f_{k-1}(x)\dots f_1(x)
$$
and if $d_k$ is a root of $P_k(x)$, then $d_k$ is a root of $f_{k-1}(x)\dots f_1(x)$.

Hence by Proposition \ref{prop:P=LR},
$P_k(d_k)d_k(P_k(d_k))^{-1}$ is a root of $P(x)$  for $k=1,\dots,n$.

By Lemma \ref{lm:linear factors}, $P(x)$ has no other roots.

\qed


In the notation of Theorem \ref{th:roots_division}, we have the following


\begin{corol} \label{cor:}
Let $\cD$ be a division ring with center $F$, $d_1,d_2 \in \cD$ such that the conjugacy classes $[d_1]$ and $[d_2]$ are distinct.
Assume that $d_1$ is algebraic over $F$.
Let also $f(x)$ be the minimal polynomial of $d_1$ and $S(x) \in \cD[x]$ such that $f(x)=S(x)(x-d_1)$.
Then
$$
(x-d_2)(x-d_1) = (x - d)(x-S(d_2)d_2(S(d_2))^{-1}),
$$
where $d= (d_1-S(d_2)d_2(S(d_2))^{-1}) d_1 (d_1-S(d_2)d_2(S(d_2))^{-1})^{-1}$.
\end{corol}

\noindent {\it {Proof}}.
Let $P(x) := (x-d_2)(x-d_1)$. By Theorem \ref{th:roots_division}, $d_3 := S(d_2)d_2(S(d_2))^{-1}$ is a root of $P(x)$. Then $x-d_3$ is a right divisor
of $P(x)$ and $P(x)=(x-d)(x-d_3)$ for some $d\in \cD$. Since $d_1$ is a root of $P(x)$ and $d_1\ne d_3$, then by Proposition \ref{prop:P=LR},
$$
d= (d_1-d_3)d_1(d_1-d_3)^{-1} = (d_1-S(d_2)d_2(S(d_2))^{-1}) d_1 (d_1-S(d_2)d_2(S(d_2))^{-1})^{-1}.
$$
\qed


\begin{remark}
The formula from the previous corollary allows to
change the order of factors in products of monic linear polynomials. This formula generalizes formulas for Hamilton's quaternion algebra from
\cite[Lm. 1]{SeSi01} (see also \cite[Th.7]{FaMiSeSo17}).
\end{remark}


\begin{example} \label{ex:field}
Let $F$ be a field, $\char(F)\ne 2$. Let also $\cQ$ be a quaternion division algebra over $F$. Assume that $d_1,d_2 \in \cQ$, $[d_1]\ne [d_2]$.
If $d_1 \not \in F$, then the minimal polynomial of $d_1$ is $(x-\overline{d}_1)(x-d_1)$,
where $\overline{d}_1$ is the conjugate of $d_1$.  Hence in the notation of Corollary \ref{cor:},
$S(x)=x-\overline{d}_1$. Then
$$
S(d_2)d_2(S(d_2))^{-1}= (d_2-\overline{d}_1)d_2(d_2-\overline{d}_1)^{-1}.
$$
Simple computations show that
$$
(d_1-S(d_2)d_2(S(d_2))^{-1}) d_1 (d_1-S(d_2)d_2(S(d_2))^{-1})^{-1} = (d_2-\overline{d}_1)d_1(d_2-\overline{d}_1)^{-1}.
$$
Thus
$$
(x-d_2)(x-d_1)= (x-hd_1h^{-1})(x-hd_2h^{-1}),
$$
where $h=d_2-\overline{d}_1$ (compare with the formula from \cite[Lm. 1]{SeSi01}).
\end{example}




\begin{thebibliography}{99}

\frenchspacing
\bibitem{Ab09} M. Abrate. {\em Quadratic formulas for generalized quaternions}, J.
Algebra Appl. {\bf 8} (2009), no. 3, 289-306. MR2535990 

\bibitem{Ch15} A. Chapman. {\em Quaternion quadratic equations in characteristic
2}, J. Algebra Appl. {\bf 14} (2015), no. 3, 1550033.  MR3275570

\bibitem{Ch20} A. Chapman. {\em Polynomial equations over octonion algebras},
J. Algebra Appl. {\bf 19} (2020), no. 6, 2050102. MR4120079

\bibitem{ChMa16}
A. Chapman, C. Machen. {\em Standard polynomial equations over division algebras},
Adv. Appl. Clifford Algebr. {\bf 27} (2017), no.2, 1065-1072. MR3651503 

\bibitem{FaMiSeSo17}
M. I. Falc\~{a}o, F. Miranda, R. Severino, M. J. Soares, {\em Mathematica Tools for Quaternionic Polynomials},
Computational science and its applications.  ICCSA 2017. Part II, 394-408.
Lecture Notes in Comput. Sci., 10405
Springer, Cham, 2017. MR3696911

\bibitem{GoMo65} B. Gordon, T.S. Motzkin. {\em On the zeros of polynomials over division rings}, Trans. Amer.
Math. Soc., {\bf 116} (1965), 218-226. MR0195853

\bibitem{HuSo02}
L. Huang,  W. So. {\em Quadratic formulas for quaternions}, Appl. Math. Lett. {\bf 15} (2002), no. 5, 533-540. MR 1889501

\bibitem{JaOp10} D. Janovsk\'{a}, G. Opfer. {\em A note on the computation of all zeros of simple quaternionic polynomials},
SIAM J. Numer. Anal. {\bf 48} (2010), no. 1, 244-256. MR2608368 

\bibitem{La91}
T.Y. Lam, {\em A first course in noncommutative rings}, Graduate Texts in Mathematics 131, Springer-Verlag, New York, 1991.
MR1125071


\bibitem{SePeVi01} R. Ser\^{o}dio, E. Pereira, J. Vit\'{o}ria, {\em Computing the zeros of
quaternion polynomials}, Comput. Math. Appl. {\bf 42} (2001), no. 8-9, 1229-1237. MR 1851239

\bibitem{SeSi01} R. Ser\^{o}dio, L.-S. Siu. {\em Zeros of quaternion polynomials},
Appl. Math. Lett. {\bf 14} (2001),  no. 2, 237-239. MR1808271


 \end{thebibliography}
\end{document}